\documentclass[12pt, reqno]{amsart}
\usepackage{amssymb}
\usepackage{eucal}
\usepackage{amsmath}
\usepackage{amscd}
\usepackage[dvips]{color}
\usepackage{multicol}
\usepackage[all]{xy}           
\usepackage{graphicx}
\usepackage{color}
\usepackage{colordvi}
\usepackage{xspace, ulem, lipsum}
\usepackage{tikz}
\usepackage{ifpdf}
\ifpdf
 \usepackage[colorlinks,final,backref=page,hyperindex]{hyperref}
\else
 \usepackage[colorlinks,final,backref=page,hyperindex,hypertex]{hyperref}
\fi

\begin{document}
\newcommand {\emptycomment}[1]{} 

\baselineskip=14pt
\newcommand{\nc}{\newcommand}
\newcommand{\delete}[1]{}
\nc{\mfootnote}[1]{\footnote{#1}} 
\nc{\todo}[1]{\tred{To do:} #1}

\delete{
\nc{\mlabel}[1]{\label{#1}}  
\nc{\mcite}[1]{\cite{#1}}  
\nc{\mref}[1]{\ref{#1}}  
\nc{\meqref}[1]{\ref{#1}} 
\nc{\mbibitem}[1]{\bibitem{#1}} 
}

\nc{\mlabel}[1]{\label{#1}  
{\hfill \hspace{1cm}{\bf{{\ }\hfill(#1)}}}}
\nc{\mcite}[1]{\cite{#1}{{\bf{{\ }(#1)}}}}  
\nc{\mref}[1]{\ref{#1}{{\bf{{\ }(#1)}}}}  
\nc{\meqref}[1]{\eqref{#1}{{\bf{{\ }(#1)}}}} 
\nc{\mbibitem}[1]{\bibitem[\bf #1]{#1}} 

\newtheorem{theorem}{Theorem}[section]
\newtheorem{lemma}[theorem]{Lemma}
\newtheorem{corrollary}[theorem]{Corollary}
\newtheorem{proposition}[theorem]{Proposition}
\theoremstyle{definition}
\newtheorem{definition}[theorem]{Definition}
\newtheorem{example}[theorem]{Example}
\newtheorem{remark}[theorem]{Remark}
\newtheorem{pdef}[theorem]{Proposition-Definition}
\newtheorem{condition}[theorem]{Condition}

\renewcommand{\labelenumi}{{\rm(\alph{enumi})}}
\renewcommand{\theenumi}{\alph{enumi}}

\nc{\tred}[1]{\textcolor{red}{#1}}
\nc{\tblue}[1]{\textcolor{blue}{#1}}
\nc{\tgreen}[1]{\textcolor{green}{#1}}
\nc{\tpurple}[1]{\textcolor{purple}{#1}}
\nc{\btred}[1]{\textcolor{red}{\bf #1}}
\nc{\btblue}[1]{\textcolor{blue}{\bf #1}}
\nc{\btgreen}[1]{\textcolor{green}{\bf #1}}
\nc{\btpurple}[1]{\textcolor{purple}{\bf #1}}

\nc{\ld}[1]{\textcolor{purple}{Landry:#1}}
\nc{\hd}[1]{\textcolor{blue}{Damien:#1}}
\nc{\li}[1]{\textcolor{yellow}{#1}}
\nc{\lir}[1]{\textcolor{blue}{Li:#1}}


\nc{\twovec}[2]{\left(\begin{array}{c} #1 \\ #2\end{array} \right )}
\nc{\threevec}[3]{\left(\begin{array}{c} #1 \\ #2 \\ #3 \end{array}\right )}
\nc{\twomatrix}[4]{\left(\begin{array}{cc} #1 & #2\\ #3 & #4 \end{array} \right)}
\nc{\threematrix}[9]{{\left(\begin{matrix} #1 & #2 & #3\\ #4 & #5 & #6 \\ #7 & #8 & #9 \end{matrix} \right)}}
\nc{\twodet}[4]{\left|\begin{array}{cc} #1 & #2\\ #3 & #4 \end{array} \right|}

\nc{\rk}{\mathrm{r}}
\newcommand{\g}{\mathfrak g}
\newcommand{\h}{\mathfrak h}
\newcommand{\pf}{\noindent{$Proof$.}\ }
\newcommand{\frkg}{\mathfrak g}
\newcommand{\frkh}{\mathfrak h}
\newcommand{\Id}{\rm{Id}}
\newcommand{\gl}{\mathfrak {gl}}
\newcommand{\ad}{\mathrm{ad}}
\newcommand{\add}{\frka\frkd}
\newcommand{\frka}{\mathfrak a}
\newcommand{\frkb}{\mathfrak b}
\newcommand{\frkc}{\mathfrak c}
\newcommand{\frkd}{\mathfrak d}
\newcommand {\comment}[1]{{\marginpar{*}\scriptsize\textbf{Comments:} #1}}

\nc{\tforall}{\text{ for all }}

\nc{\svec}[2]{{\tiny\left(\begin{matrix}#1\\
#2\end{matrix}\right)\,}}  
\nc{\ssvec}[2]{{\tiny\left(\begin{matrix}#1\\
#2\end{matrix}\right)\,}} 

\nc{\typeI}{local cocycle $3$-Lie bialgebra\xspace}
\nc{\typeIs}{local cocycle $3$-Lie bialgebras\xspace}
\nc{\typeII}{double construction $3$-Lie bialgebra\xspace}
\nc{\typeIIs}{double construction $3$-Lie bialgebras\xspace}

\nc{\bia}{{$\mathcal{P}$-bimodule ${\bf k}$-algebra}\xspace}
\nc{\bias}{{$\mathcal{P}$-bimodule ${\bf k}$-algebras}\xspace}

\nc{\rmi}{{\mathrm{I}}}
\nc{\rmii}{{\mathrm{II}}}
\nc{\rmiii}{{\mathrm{III}}}
\nc{\pr}{{\mathrm{pr}}}
\newcommand{\huaA}{\mathcal{A}}

\nc{\OT}{constant $\theta$-}
\nc{\T}{$\theta$-}
\nc{\IT}{inverse $\theta$-}

\nc{\pll}{\beta}
\nc{\plc}{\epsilon}

\nc{\ass}{{\mathit{Ass}}}
\nc{\lie}{{\mathit{Lie}}}
\nc{\comm}{{\mathit{Comm}}}
\nc{\dend}{{\mathit{Dend}}}
\nc{\zinb}{{\mathit{Zinb}}}
\nc{\tdend}{{\mathit{TDend}}}
\nc{\prelie}{{\mathit{preLie}}}
\nc{\postlie}{{\mathit{PostLie}}}
\nc{\quado}{{\mathit{Quad}}}
\nc{\octo}{{\mathit{Octo}}}
\nc{\ldend}{{\mathit{ldend}}}
\nc{\lquad}{{\mathit{LQuad}}}

 \nc{\adec}{\check{;}} \nc{\aop}{\alpha}
\nc{\dftimes}{\widetilde{\otimes}} \nc{\dfl}{\succ} \nc{\dfr}{\prec}
\nc{\dfc}{\circ} \nc{\dfb}{\bullet} \nc{\dft}{\star}
\nc{\dfcf}{{\mathbf k}} \nc{\apr}{\ast} \nc{\spr}{\cdot}
\nc{\twopr}{\circ} \nc{\tspr}{\star} \nc{\sempr}{\ast}
\nc{\disp}[1]{\displaystyle{#1}}
\nc{\bin}[2]{ (_{\stackrel{\scs{#1}}{\scs{#2}}})}  
\nc{\binc}[2]{ \left (\!\! \begin{array}{c} \scs{#1}\\
    \scs{#2} \end{array}\!\! \right )}  
\nc{\bincc}[2]{  \left ( {\scs{#1} \atop
    \vspace{-.5cm}\scs{#2}} \right )}  
\nc{\sarray}[2]{\begin{array}{c}#1 \vspace{.1cm}\\ \hline
    \vspace{-.35cm} \\ #2 \end{array}}
\nc{\bs}{\bar{S}} \nc{\dcup}{\stackrel{\bullet}{\cup}}
\nc{\dbigcup}{\stackrel{\bullet}{\bigcup}} \nc{\etree}{\big |}
\nc{\la}{\longrightarrow} \nc{\fe}{\'{e}} \nc{\rar}{\rightarrow}
\nc{\dar}{\downarrow} \nc{\dap}[1]{\downarrow
\rlap{$\scriptstyle{#1}$}} \nc{\uap}[1]{\uparrow
\rlap{$\scriptstyle{#1}$}} \nc{\defeq}{\stackrel{\rm def}{=}}
\nc{\dis}[1]{\displaystyle{#1}} \nc{\dotcup}{\,
\displaystyle{\bigcup^\bullet}\ } \nc{\sdotcup}{\tiny{
\displaystyle{\bigcup^\bullet}\ }} \nc{\hcm}{\ \hat{,}\ }
\nc{\hcirc}{\hat{\circ}} \nc{\hts}{\hat{\shpr}}
\nc{\lts}{\stackrel{\leftarrow}{\shpr}}
\nc{\rts}{\stackrel{\rightarrow}{\shpr}} \nc{\lleft}{[}
\nc{\lright}{]} \nc{\uni}[1]{\tilde{#1}} \nc{\wor}[1]{\check{#1}}
\nc{\free}[1]{\bar{#1}} \nc{\den}[1]{\check{#1}} \nc{\lrpa}{\wr}
\nc{\curlyl}{\left \{ \begin{array}{c} {} \\ {} \end{array}
    \right .  \!\!\!\!\!\!\!}
\nc{\curlyr}{ \!\!\!\!\!\!\!
    \left . \begin{array}{c} {} \\ {} \end{array}
    \right \} }
\nc{\leaf}{\ell}       
\nc{\longmid}{\left | \begin{array}{c} {} \\ {} \end{array}
    \right . \!\!\!\!\!\!\!}
\nc{\ot}{\otimes} \nc{\sot}{{\scriptstyle{\ot}}}
\nc{\otm}{\overline{\ot}}
\nc{\ora}[1]{\stackrel{#1}{\rar}}
\nc{\ola}[1]{\stackrel{#1}{\la}}
\nc{\pltree}{\calt^\pl}
\nc{\epltree}{\calt^{\pl,\NC}}
\nc{\rbpltree}{\calt^r}
\nc{\scs}[1]{\scriptstyle{#1}} \nc{\mrm}[1]{{\rm #1}}
\nc{\dirlim}{\displaystyle{\lim_{\longrightarrow}}\,}
\nc{\invlim}{\displaystyle{\lim_{\longleftarrow}}\,}
\nc{\mvp}{\vspace{0.5cm}} \nc{\svp}{\vspace{2cm}}
\nc{\vp}{\vspace{8cm}} \nc{\proofbegin}{\noindent{\bf Proof: }}
\nc{\proofend}{$\blacksquare$ \vspace{0.5cm}}
\nc{\freerbpl}{{F^{\mathrm RBPL}}}
\nc{\sha}{{\mbox{\cyr X}}}  
\nc{\ncsha}{{\mbox{\cyr X}^{\mathrm NC}}} \nc{\ncshao}{{\mbox{\cyr
X}^{\mathrm NC,\,0}}}
\nc{\shpr}{\diamond}    
\nc{\shprm}{\overline{\diamond}}    
\nc{\shpro}{\diamond^0}    
\nc{\shprr}{\diamond^r}     
\nc{\shpra}{\overline{\diamond}^r}
\nc{\shpru}{\check{\diamond}} \nc{\catpr}{\diamond_l}
\nc{\rcatpr}{\diamond_r} \nc{\lapr}{\diamond_a}
\nc{\sqcupm}{\ot}
\nc{\lepr}{\diamond_e} \nc{\vep}{\varepsilon} \nc{\labs}{\mid\!}
\nc{\rabs}{\!\mid} \nc{\hsha}{\widehat{\sha}}
\nc{\lsha}{\stackrel{\leftarrow}{\sha}}
\nc{\rsha}{\stackrel{\rightarrow}{\sha}} \nc{\lc}{\lfloor}
\nc{\rc}{\rfloor}
\nc{\tpr}{\sqcup}
\nc{\nctpr}{\vee}
\nc{\plpr}{\star}
\nc{\rbplpr}{\bar{\plpr}}
\nc{\sqmon}[1]{\langle #1\rangle}
\nc{\forest}{\calf}
\nc{\altx}{\Lambda_X} \nc{\vecT}{\vec{T}} \nc{\onetree}{\bullet}
\nc{\Ao}{\check{A}}
\nc{\seta}{\underline{\Ao}}
\nc{\deltaa}{\overline{\delta}}
\nc{\trho}{\tilde{\rho}}

\nc{\rpr}{\circ}
\nc{\dpr}{{\tiny\diamond}}
\nc{\rprpm}{{\rpr}}

\nc{\mmbox}[1]{\mbox{\ #1\ }} \nc{\ann}{\mrm{ann}}
\nc{\Aut}{\mrm{Aut}} \nc{\can}{\mrm{can}}
\nc{\twoalg}{{two-sided algebra}\xspace}
\nc{\colim}{\mrm{colim}}
\nc{\Cont}{\mrm{Cont}} \nc{\rchar}{\mrm{char}}
\nc{\cok}{\mrm{coker}} \nc{\dtf}{{R-{\rm tf}}} \nc{\dtor}{{R-{\rm
tor}}}
\renewcommand{\det}{\mrm{det}}
\nc{\depth}{{\mrm d}}
\nc{\Div}{{\mrm Div}} \nc{\End}{\mrm{End}} \nc{\Ext}{\mrm{Ext}}
\nc{\Fil}{\mrm{Fil}} \nc{\Frob}{\mrm{Frob}} \nc{\Gal}{\mrm{Gal}}
\nc{\GL}{\mrm{GL}} \nc{\Hom}{\mrm{Hom}} \nc{\hsr}{\mrm{H}}
\nc{\hpol}{\mrm{HP}} \nc{\id}{\mrm{id}} \nc{\im}{\mrm{im}}
\nc{\incl}{\mrm{incl}} \nc{\length}{\mrm{length}}
\nc{\LR}{\mrm{LR}} \nc{\mchar}{\rm char} \nc{\NC}{\mrm{NC}}
\nc{\mpart}{\mrm{part}} \nc{\pl}{\mrm{PL}}
\nc{\ql}{{\QQ_\ell}} \nc{\qp}{{\QQ_p}}
\nc{\rank}{\mrm{rank}} \nc{\rba}{\rm{RBA }} \nc{\rbas}{\rm{RBAs }}
\nc{\rbpl}{\mrm{RBPL}}
\nc{\rbw}{\rm{RBW }} \nc{\rbws}{\rm{RBWs }} \nc{\rcot}{\mrm{cot}}
\nc{\rest}{\rm{controlled}\xspace}
\nc{\rdef}{\mrm{def}} \nc{\rdiv}{{\rm div}} \nc{\rtf}{{\rm tf}}
\nc{\rtor}{{\rm tor}} \nc{\res}{\mrm{res}} \nc{\SL}{\mrm{SL}}
\nc{\Spec}{\mrm{Spec}} \nc{\tor}{\mrm{tor}} \nc{\Tr}{\mrm{Tr}}
\nc{\mtr}{\mrm{sk}}

\nc{\ab}{\mathbf{Ab}} \nc{\Alg}{\mathbf{Alg}}
\nc{\Algo}{\mathbf{Alg}^0} \nc{\Bax}{\mathbf{Bax}}
\nc{\Baxo}{\mathbf{Bax}^0} \nc{\RB}{\mathbf{RB}}
\nc{\RBo}{\mathbf{RB}^0} \nc{\BRB}{\mathbf{RB}}
\nc{\Dend}{\mathbf{DD}} \nc{\bfk}{{\bf k}} \nc{\bfone}{{\bf 1}}
\nc{\base}[1]{{a_{#1}}} \nc{\detail}{\marginpar{\bf More detail}
    \noindent{\bf Need more detail!}
    \svp}
\nc{\Diff}{\mathbf{Diff}} \nc{\gap}{\marginpar{\bf
Incomplete}\noindent{\bf Incomplete!!}
    \svp}
\nc{\FMod}{\mathbf{FMod}} \nc{\mset}{\mathbf{MSet}}
\nc{\rb}{\mathrm{RB}} \nc{\Int}{\mathbf{Int}}
\nc{\Mon}{\mathbf{Mon}}
\nc{\remarks}{\noindent{\bf Remarks: }}
\nc{\OS}{\mathbf{OS}} 
\nc{\Rep}{\mathbf{Rep}}
\nc{\Rings}{\mathbf{Rings}} \nc{\Sets}{\mathbf{Sets}}
\nc{\DT}{\mathbf{DT}}

\nc{\BA}{{\mathbb A}} \nc{\CC}{{\mathbb C}} \nc{\DD}{{\mathbb D}}
\nc{\EE}{{\mathbb E}} \nc{\FF}{{\mathbb F}} \nc{\GG}{{\mathbb G}}
\nc{\HH}{{\mathbb H}} \nc{\LL}{{\mathbb L}} \nc{\NN}{{\mathbb N}}
\nc{\QQ}{{\mathbb Q}} \nc{\RR}{{\mathbb R}} \nc{\BS}{{\mathbb{S}}} \nc{\TT}{{\mathbb T}}
\nc{\VV}{{\mathbb V}} \nc{\ZZ}{{\mathbb Z}}


\nc{\calao}{{\mathcal A}} \nc{\cala}{{\mathcal A}}
\nc{\calc}{{\mathcal C}} \nc{\cald}{{\mathcal D}}
\nc{\cale}{{\mathcal E}} \nc{\calf}{{\mathcal F}}
\nc{\calfr}{{{\mathcal F}^{\,r}}} \nc{\calfo}{{\mathcal F}^0}
\nc{\calfro}{{\mathcal F}^{\,r,0}} \nc{\oF}{\overline{F}}
\nc{\calg}{{\mathcal G}} \nc{\calh}{{\mathcal H}}
\nc{\cali}{{\mathcal I}} \nc{\calj}{{\mathcal J}}
\nc{\call}{{\mathcal L}} \nc{\calm}{{\mathcal M}}
\nc{\caln}{{\mathcal N}} \nc{\calo}{{\mathcal O}}
\nc{\calp}{{\mathcal P}} \nc{\calq}{{\mathcal Q}} \nc{\calr}{{\mathcal R}}
\nc{\calt}{{\mathcal T}} \nc{\caltr}{{\mathcal T}^{\,r}}
\nc{\calu}{{\mathcal U}} \nc{\calv}{{\mathcal V}}
\nc{\calw}{{\mathcal W}} \nc{\calx}{{\mathcal X}}
\nc{\CA}{\mathcal{A}}

\nc{\fraka}{{\mathfrak a}} \nc{\frakB}{{\mathfrak B}}
\nc{\frakb}{{\mathfrak b}} \nc{\frakd}{{\mathfrak d}}
\nc{\oD}{\overline{D}}
\nc{\frakF}{{\mathfrak F}} \nc{\frakg}{{\mathfrak g}}
\nc{\frakm}{{\mathfrak m}} \nc{\frakM}{{\mathfrak M}}
\nc{\frakMo}{{\mathfrak M}^0} \nc{\frakp}{{\mathfrak p}}
\nc{\frakS}{{\mathfrak S}} \nc{\frakSo}{{\mathfrak S}^0}
\nc{\fraks}{{\mathfrak s}} \nc{\os}{\overline{\fraks}}
\nc{\frakT}{{\mathfrak T}}
\nc{\oT}{\overline{T}}
\nc{\frakX}{{\mathfrak X}} \nc{\frakXo}{{\mathfrak X}^0}
\nc{\frakx}{{\mathbf x}}
\nc{\frakTx}{\frakT}      
\nc{\frakTa}{\frakT^a}        
\nc{\frakTxo}{\frakTx^0}   
\nc{\caltao}{\calt^{a,0}}   
\nc{\ox}{\overline{\frakx}} \nc{\fraky}{{\mathfrak y}}
\nc{\frakz}{{\mathfrak z}} \nc{\oX}{\overline{X}}

\font\cyr=wncyr10

\nc{\al}{\alpha}
\nc{\lam}{\lambda}
\nc{\lr}{\longrightarrow}
\newcommand{\K}{\mathbb {K}}
\newcommand{\A}{\rm A}


\title[On Rota-Baxter and Nijenhuis anti-flexible algebras]{On Rota-Baxter and Nijenhuis anti-flexible algebras}

\author[Mafoya Landry Dassoundo]{Mafoya Landry Dassoundo}
\address[]{Department of Mathematics, University of Man,  BP 20 Man,  Côte d'Ivoire (Ivory Coast) \&
	University of Abomey-Calavi, International Chair 
	in Mathematical Physics and its Applications, 
	ICMPA-UNESCO Chair, 072 BP 50, Cotonou, Rep. of Benin}
\email{landry.dassoundo@univ-man.edu.ci \\ mafoya.dassoundo@cipma.uac.bj}

\begin{abstract}
	We define and derive basic properties of the
	notion of Rota-Baxter operator on anti-flexible algebra.
	Starting from a Rota-Baxter operator on an anti-flexible algebra, 
	we construct pre-anti-flexible algebra structure and associated 
	left(right)-symmetric algebra as well.
	The notion of $\mathcal{O}$-operator on anti-flexible 
	algebra is recalled  and used to build left(right)-symmetric
	algebra as well as  related properties. Furthermore, we introduce
	Nijenhuis anti-flexible algebra and 
	derive associated properties. Nijenhuis operator on
	anti-flexible algebra is used to build pre-anti-flexible 
	algebra structure and related left(right)-symmetric algebra.
\end{abstract}
\maketitle
\section{Introduction and Preliminaries}
First of all, we recall the definition and useful properties of ﬂexible algebra.
\begin{definition}\label{def.-flexible1}
	A flexible algebra is a $\mathbb{K}$-vector 
	space $A$ equipped with a bilinear product $\ast: A\times A\rightarrow A$
	satisfying  
	\begin{eqnarray}\label{eq:flexible1}
		[x,y,x]:=	(x\ast y)\ast x-x\ast (y\ast x)=0, \quad  \forall x,y\in A.
	\end{eqnarray}
\end{definition}
Substituting $x$ in Eq.~\eqref{eq:flexible1} by $x+z$, we then derive the following equation  
\begin{eqnarray*}
	[x,y,z]+[z, y, x]=0, \forall x,y,z\in A.
\end{eqnarray*}
Clearly, as the characteristic of the considered 
field $\mathbb{K}$ in this paper is
different to 2, the following definition is then equivalent
to Definition~\ref{def.-flexible1}.
\begin{definition}
	A flexible algebra is a $\mathbb{K}$-vector 
	space $A$ equipped with a bilinear product $\ast: A\times A\rightarrow A$
	satisfying 
	\begin{eqnarray}\label{eq:flexible}
		[x,y,z]+[z, y, x]=0, \forall x,y,z\in A.
	\end{eqnarray}
\end{definition}
Introduced as a natural generalization of associative algebras~\cite{Oehmke}, 
flexible algebras are also Cayley-Dickson algebras~\cite{Schafer}.
A generalization of flexible algebras leads to 
several classes of non-associative algebras
among them we have the so call
anti-flexible algebras~\cite{Rodabaugh_3} also know as 
center-symmetric algebras and investigated  \cite{Hounkonnou_D_CSA}, 
the flexible Lie-admissible algebras which are studied
in \cite{Benkart_O} and the references therein.
\begin{definition}\cite[Definition 1.1]{Rodabaugh_3}
	An anti-flexible algebra is a vector space  $A$ equipped with a linear product
	$\cdot: A\otimes A\rightarrow A, x\otimes y\mapsto x\cdot y,$ satisfying  
	\begin{eqnarray}
		(x\cdot y)\cdot z-x\cdot (y\cdot z)=
		(z\cdot y)\cdot x-z\cdot (y\cdot x), \forall x,y,z\in A,
	\end{eqnarray}	
	or equivalently
	\begin{eqnarray}
		[x,y,z]=
		[z,y,x], \forall x,y,z\in A.
	\end{eqnarray}		
\end{definition}
The study of non-flexible
algebras which are also  quasi-associative algebras
revealed their connections with Jordan-admissible and Lie-admissible  algebras \cite{Kosier} and the references therein.
Investigation on algebraic properties of anti-flexible  
algebras gives rise to
interesting results. For instance, 
the complete classification of simple 
anti-flexible algebras on a field of 
characteristic different to $2$ and $3,$ 
as well as some characterization of 
semi-simple anti-flexible algebras
\cite{Bhandari, Rodabaugh_1, Rodabaugh_2, Rodabaugh}.
In addition,  it is important to notice that the adjoint functor between the category of 
Rota-Baxter algebras and the categories of dendriform
dialgebras and trialgebras were investigated  and 
related properties were derived in \cite{Ebrahimi-Fard_G}.
\begin{proposition}\label{pro:opp}
	Consider an anti-flexible algebra $(A, \mu),$ where $\mu: A\otimes A\rightarrow A, x\otimes y\mapsto \mu(x\otimes y),$ and 
	$\sigma: A\otimes A\rightarrow A\otimes A, x\otimes y\mapsto y\otimes x.$
	We have:
	\begin{itemize}
		\item[$\bullet$] 
		The opposite algebra $A^{op}:=(A, \mu\circ \sigma)$ is also 
		an anti-flexible algebra.
		\item[$\bullet$] 
		The algebra $(A, \mu-\mu\circ\sigma)$ is a Lie algebra, where  
		$\forall x,y\in A,$ \\$(\mu-\mu\circ\sigma)(x\otimes y)=\mu(x\otimes y)-\mu(y\otimes x).$
	\end{itemize}
\end{proposition}
\begin{example}
	For a given vector space $V$	equipped with a scalar product 
	$<,>:A\times A\rightarrow \mathbb{K},$ and with a linear form
	$L: A\rightarrow \mathbb{K}$ with $L(c)=0,$ where $c$ is   a fix element 
	$c\in A.$  
	There is an anti-flexible algebra structure on $V$ given by, 
	\begin{eqnarray*}
		x\cdot y= <x,c><y,c>c+L(x)y, \forall x,y\in V.
	\end{eqnarray*}
	Since the scalar product commute on $V\times V,$ then in view of Proposition~\ref{pro:opp}, the following 
	product given by,  
	\begin{eqnarray*}
		x\star y= <x,c><y,c>c+L(y)x, \forall x,y\in V
	\end{eqnarray*}
	together with $V$ define an anti-flexible algebra.
	In particular, if $c_1\in \{x\in V, <x,c>=0\},$ then the 
	following products given by  $\forall x,y\in V,$ 
	\begin{subequations}
		\begin{eqnarray*}
			x\cdot_1 y= <x,c><y,c>c+<x, c_1>y;
		\end{eqnarray*}
		\begin{eqnarray*}
			x\cdot_2 y= <x,c><y,c>c+<y,c_1>x
		\end{eqnarray*}
	\end{subequations}
	define anti-flexible algebra on $V.$
\end{example}
\begin{definition}\label{Def1}\cite[Definition 5.2]{Dassoundo_B_H}
	A pre-anti-flexible algebra is a vector space 
	$A$ equipped with two linear maps 
	$\prec, \succ: A\times A\rightarrow A$ 
	satisfying the following identities, 
	$\forall x,y, z\in A,$
	\begin{subequations}
		\begin{eqnarray}\label{eq:pre-anti-flexible1}
			(x\succ y)\prec z-x\succ (y\prec z)=
			(z\succ y)\prec x- z\succ (y\prec x), 
		\end{eqnarray}
		\begin{eqnarray}\label{eq:pre-anti-flexible2}
			(x\succ y+x\prec y)\succ z-x\succ (y\succ z)=\cr
			(z\prec y)\prec x-z\prec (y\prec x+y\succ x).
		\end{eqnarray}
	\end{subequations} 
\end{definition}
\begin{remark}\label{rmk:one}
	By considering the Definition~\ref{Def1}, we have:
	\begin{enumerate}
		\item \label{one}
		It follows from Eq.~\eqref{eq:pre-anti-flexible1} and 
		Eq.~\eqref{eq:pre-anti-flexible2}
		that the bilinear product $\ast: A\rightarrow A$ defined by 
		$x\ast y= x\succ y+x\prec y$ turns $(A, \ast)$ into an anti-flexible algebra. 
		\item \label{two}
		By vanishing the both side of the equal sign in 
		Eq.~\eqref{eq:pre-anti-flexible1} and 
		Eq.~\eqref{eq:pre-anti-flexible2}, the pre-anti-flexible algebra $(A, \prec, \succ)$
		turns into a dendriform algebra which is introduced by Ebrahimi-Fard \& Guo in \cite{Ebrahimi-Fard_G}.
		Clearly, pre-anti-flexible algebras 
		are a generalization of dendriform algebras.
	\end{enumerate}
\end{remark}
Furthermore,  a general operadic definition for the notion of 
splitting the products of algebraic structures is provided and it has been 
proved that some Manin products
of operads are equivalent and also closed to  Rota-Baxter operators \cite{BBGN}.
\begin{definition}\cite{Vinberg}
	A couple $(A, \cdot)$, where $A$ is a vector space and $\cdot: A\times A\rightarrow A$
	is a linear map on $A$, is called 
	left-symmetric (resp. right-symmetric) algebra if the following condition 
	\begin{eqnarray*}
		(x\cdot y)\cdot z-x\cdot(y\cdot z)= (y\cdot x)\cdot z-y\cdot (x\cdot z), \forall x,y,z\in A,
	\end{eqnarray*}
	\begin{eqnarray*}
		(x\cdot y)\cdot z-x\cdot(y\cdot z)= (x\cdot z)\cdot y-x\cdot (z\cdot y), \forall x,y,z\in A,
	\end{eqnarray*}
	respectively, is satisfied.
\end{definition}
\begin{theorem}\cite[Theorem~2.8]{Mafoya}\label{thm:Prelie}
	Let $(A, \prec,  \succ)$ be a pre-anti-flexible algebra. We have 
	\begin{itemize}
		\item $(A, \triangleright)$ is a left-symmetric algebra, 
		where  $\forall x,y\in A,$
		\begin{eqnarray*}
			x\triangleright y= x\succ y-y\prec x.
		\end{eqnarray*}
		\item $(A, \triangleleft)$ is a right-symmetric algebra, 
		where  $\forall x, y\in A,$
		\begin{eqnarray*}
			x\triangleleft y=x\prec y-y\succ x.
		\end{eqnarray*}
	\end{itemize}
\end{theorem}
\begin{remark}
	As we highlighted in Remark~\ref{rmk:one} (\ref{two}) that pre-anti-flexible algebras
	generalize dendriform algebras, left (right) symmetric algebras 
	derived in Theorem~\ref{thm:Prelie} 
	take into account those constructed from dendriform 
	algebras in \cite[page 264]{Aguiar}. 	
\end{remark}
Besides, zero weight Rota-Baxter algebras and dendriform algebras are linked in \cite{Aguiar}  such away that zero weight Rota-Baxter algebras carry  dendriform dialgebras structures which look like  Lie algebras structures on associative algebras. 

It has already been established in \cite{Dassoundo_B_H} that anti-flexible algebras are a generalization of associative algebras, and associated Lie algebras are similar to Lie algebras related to associative algebras.

In addition, A. Das in \cite{Das} constructed a cohomology for operator $\mathcal{O}$-operator on associative algebras and introduced Nijenhuis elements associated with an
$\mathcal{O}$-operator on associative algebras. 
It is then worth to investigate Rota-Baxter anti-flexible algebras and pre-anti-flexible algebras which are a generalization of dendriform algebras, and also investigate on  Nijenhuis structures on anti-flexible algebras  and their relation with pre-anti-flexible algebras.
{
	The content of  this paper is as follow. In section~\ref{section2},  we provide  basic definitions 
	and derive basic properties of the notion of Rota-Baxter anti-flexible algebras,
	and  under some conditions, we construct a pre-anti-flexible algebra  structure related to a given 
	Rota-Baxter operator on an anti-flexible algebra 
	and  deduce its associated left(right)-symmetric algebra. 
	Furthermore, in section~\ref{section3}, the notion of bimodule of 
	anti-flexible algebras are recalled as well as  that of
	$\mathcal{O}$-operators on anti-flexible algebras. In addition, pre-anti-flexible algebras
	are constructed from a given $\mathcal{O}$-operator on anti-flexible
	and the induced left(right)-symmetric algebra are derived as well as related consequences.
	Finally in section~\ref{section4}, basic definitions and properties 
	on Nijenhuis anti-flexible algebras were introduced and used 
	to construct pre-anti-flexible algebra and derive left(right)-symmetric 
	algebra associated.
}
\section{Anti-flexible Rota-Baxter algebras}\label{section2}
In this section, we present basic definitions and properties 
of anti-flexible Rota-Baxter algebras. We construct a
pre-anti-flexible algebra structure from a given anti-flexible
Rota-Baxter algebra and derive their associated left- and
right-symmetric algebras.
\begin{definition}
	Let $(A, \cdot)$ be an anti-flexible algebra. A linear map $R_B:A\rightarrow A$ 
	is said to be a Rota-Baxter operator of weight $\lambda\in \mathbb{K}$
	if  
	\begin{eqnarray}\label{eq:RB-operator}
		R_B(a)\cdot R_B(b)=R_B(a\cdot R_B(b)+R_B(a)\cdot b)+\lambda R_B(a\cdot b), \forall a,b\in A.
	\end{eqnarray}
\end{definition}
\begin{proposition}\label{pro:RB}
	Let $(A, \cdot)$ be an anti-flexible algebra
	and $R_B:A\rightarrow A$ be a Rota-Baxter operator 
	of weight $\lambda \in \mathbb{K}$ on $A.$
	The bilinear operation  $\cdot_{_{R_B}}: A\times A\rightarrow A$
	given by 
	\begin{eqnarray}\label{eq:product-RB}
		a\cdot_{_{R_B}} b=a\cdot R_B(b)+R_B(a)\cdot b+\lambda a\cdot b, 
		\quad \forall a,b \in A,
	\end{eqnarray}
	is such that $(A, \cdot_{_{R_B}})$ is 
	an anti-flexible algebra.
\end{proposition}
\begin{proof}
	Let $a,b,c\in A$. 
	
	The associator $[a,b,c]_{_{R_B}}:=(a\cdot_{_{R_B}}b)\cdot_{_{R_B}} c-a\cdot_{_{R_B}} (b\cdot_{_{R_B}}c )$ can be expressed as follow:
	\begin{eqnarray*}
		[a,b,c]_{_{R_B}}&=&
		(a\cdot R_B(b))\cdot R_B(c) +
		R_B(a\cdot R_B(b))\cdot c+
		\lambda(a\cdot R_B(b))\cdot c\cr 
		&&+
		(R_B(a)\cdot b)\cdot R_B(c)+
		R_B(R_B(a)\cdot b)\cdot c\cr&&+
		\lambda (R_B(a)\cdot b)\cdot c+
		\lambda(a\cdot b)\cdot R_B(c)+
		\lambda R_B(a\cdot b)\cdot c\cr&&+
		\lambda^2 (a\cdot b)\cdot c
		-a\cdot R_B(b\cdot R_B(c))
		-R_B(a)\cdot(b\cdot R_B(c))\cr&&
		-\lambda a\cdot(b\cdot R_B(c))
		-a\cdot R_B(R_B(b)\cdot c)\cr&&
		-R_B(a)\cdot(R_B(b)\cdot c)
		-\lambda a\cdot (R_B(b)\cdot c) 
		-\lambda a\cdot R_B(b\cdot c)\cr &&
		-\lambda R_B(a)\cdot (b\cdot c)
		-\lambda^2 a\cdot(b\cdot c) \cr
		[a,b,c]_{_{R_B}} &=&
		\{(R_B(a)\cdot b)\cdot R_B(c)-a\cdot R_B(b\cdot R_B(c))
		\cr&&-a\cdot R_B(R_B(b)\cdot c)-\lambda a\cdot (R_B(b)\cdot c)  \}\cr &&+
		\{R_B(a\cdot R_B(b))\cdot c+R_B(R_B(a)\cdot b)\cdot c\cr&&+\lambda R_B(a\cdot b)\cdot c -R_B(a)\cdot(R_B(b)\cdot c) \}\cr &&+
		\{(R_B(a)\cdot b)\cdot R_B(c)-R_B(a)\cdot(b\cdot R_B(c)) \}\cr&&+
		\lambda^2\{(a\cdot b)\cdot c -a\cdot(b\cdot c)\}\cr &&+
		\lambda\{(a\cdot R_B(b))\cdot c-a\cdot (R_B(b)\cdot c) \}\cr&&+
		\lambda\{(R_B(a)\cdot b)\cdot c-R_B(a)\cdot (b\cdot c) \}\cr &&+
		\lambda\{a\cdot b)\cdot R_B(c)-a\cdot(b\cdot R_B(c))\}.
	\end{eqnarray*}
	Using  Eq.~\eqref{eq:RB-operator}, the associator 
	$[a,b,c]_{_{R_B}}$
	can be rewritten into the  following form: 
	\begin{eqnarray}\label{eq:associator-product-RB}
		[a,b,c]_{_{R_B}}&=&
		[a, R_B(b), R_B(c)]+
		[R_B(a), b, R_B(c)]+
		[R_B(a), R_B(b), c]\cr&+&
		\lambda[R_B(a), b,c]+
		\lambda[a,R_B(b), c]+
		\lambda[a,b, R_B(c)]+
		\lambda^2[a,b,c].
	\end{eqnarray} 
	Since $(A, \cdot)$ is an 
	anti-flexible algebra and $R_B: A\rightarrow A$ 
	a Rota-Baxter operator of weight $\lambda$ on $A$, 
	we easily establish that 
	$[a,b,c]_{_{R_B}}=[c,b,a]_{_{R_B}}$.
	Then $(A, \cdot_{_{R_B}})$ is an 
	anti-flexible algebra.
\end{proof}
\begin{remark}
	Let $(A, \cdot)$ be an anti-flexible algebra
	and $R_B:A\rightarrow A$ be a Rota-Baxter operator 
	of weight $\lambda \in \mathbb{K}$ on $A.$
	\begin{itemize}
		\item 
		In view of Proposition~\ref{pro:RB}, 
		the Rota-Baxter operator induce
		an anti-flexible algebra homomorphism
		from $(A, \cdot_{_{R_B}})$ into 
		$(A, \cdot)$.
		\item
		According to the  form of Eq.~\eqref{eq:associator-product-RB}, 
		we deduce that the Proposition~\ref{pro:RB} remains valid
		if $(A, \cdot)$ is a pre-Lie algebra and particularly of associative algebra. 
	\end{itemize}
\end{remark}
\begin{lemma}
	Let $(A, \cdot)$ be an anti-flexible algebra 
	and $R_B:A\rightarrow A$ be a Rota-Baxter operator on $A$. 
	For any $p,q\in \mathbb{N}$ we have:
	\begin{itemize}
		\item 
		$(A, \cdot_{_{{R_B}^p}})$ is an anti-flexible algebra;
		\item 
		${R_{B}}^q$ is also a Rota-Baxter operator on the anti-flexible algebra $(A, \cdot_{_{{R_B}^p}});$
		\item 
		The anti-flexible algebras $(A, (\cdot_{_{{R_B}^p}}){_{_{{R_B}^q}})}$ and 
		$(A, \cdot_{_{{R_B}^{p+q}}})$ coincide;
		\item 
		The anti-flexible algebras $(A, \cdot_{_{{R_N}^p}})$ and $(A, \cdot_{_{{R_B}^q}})$ are compatible, 
		that is, any linear combination of 
		``$\cdot_{_{{R_B}^p}}$" and ``$\cdot_{_{{R_B}^q}}$" still induce anti-flexible algebra structure on $A$;
		\item 
		${R_{B}}^q$ is also a homomorphism from the anti-flexible algebra \newline
		$(A, \cdot_{_{{R_B}^{p+q}}})$ to $(A, \cdot_{_{{R_B}^p}})$. 
	\end{itemize}
\end{lemma}
\begin{proposition}
	Let $(A, \cdot)$ be an anti-flexible algebra 
	and  $R_B:A\rightarrow A$ be a linear map. 
	The set $\{(R_B(a), a), a\in A\}\subset R_B(A)\oplus A $ equipped with the semi-direct product given by 
	\begin{eqnarray}\label{eq:RB-subalgebraproduct}
		(R_B(a), a)\ast (R_B(b), b)&=&
		(R_B(a)\cdot R_B(b), a\cdot R_B(b)+R_B(a)\cdot b+\lambda a\cdot b), \forall a,b\in A,
	\end{eqnarray}
	is an anti-flexible algebra  if and only if
	the triple $(A, \cdot, R_B)$ is an anti-flexible 
	Rota-Baxter algebra, where $ R_B$ is a Rota-Baxter operator of weight $\lambda$ on $A$.
\end{proposition}
\begin{proof}
	Clearly, the product given in Eq.~\eqref{eq:RB-subalgebraproduct} is 
	stable on the set \newline $\{(R_B(a), a), a\in A\}$ if $R_B$ satisfies 
	Eq.~\eqref{eq:RB-operator}.
	Besides, the associator of the product defined in the equation Eq.~\eqref{eq:RB-subalgebraproduct} can be expressed as follow: 
	$\forall a,b,c\in A,[(R_B(a), a), (R_B(b),b), (R_B(c),c)]=\\
	((R_B(a), a)\ast(R_B(b), b))\ast(R_B(c), c)-
	(R_B(a), a)\ast((R_B(b), b)\ast(R_B(c), c))
	$
	can be expressed as the following form
	\begin{eqnarray}\label{eq:associator-sum}
		[(R_B(a), a), (R_B(b),b), (R_B(c),c)]=
		(R_B([a,b,c]_{_{R_B}}), [a,b,c]_{_{R_B}}).
	\end{eqnarray}
	Using the Proposition~\ref{pro:RB} and 
	Eq.~\eqref{eq:associator-sum}, we easily complete the proof.   
\end{proof}
\begin{definition}
	A morphism between two Rota-Baxter algebras \newline
	$(A, \cdot,  R_B)$ and $(A', \cdot', {R_B}')$ is 
	an anti-flexible morphism
	$\varphi:A\rightarrow A'$ such that the following
	diagram commutes,
	\begin{eqnarray}\label{eq:morphism-RB}
		\begin{array}{cccc}
			\xymatrix{	
				& A\ar[rr]^-{\varphi}\ar[d]^-{R_B}
				&& A' \ar[d]^-{{R_B}'}\cr
				& A\ar[rr]^-{\varphi}
				&& A'
			}
		\end{array}
		\; \; \mbox{  equivalently  } \; \; \; 
		\varphi \circ R_B={R_B}'\circ \varphi.
	\end{eqnarray}
\end{definition}
\begin{proposition}
	A linear map 
	$\varphi: A\rightarrow A'$
	is a morphism of the anti-flexible
	Rota-Baxter algebras
	$(A, \cdot, R_B)$ and 
	$(A', \cdot', {R_B}')$
	if and only if the set 
	\begin{eqnarray*}
		\{
		((R_B(a), a), (\varphi(R_B(a)), \varphi(a))
		)|  a\in A
		\}	\subset (R_B(A)\oplus A)\oplus ({R_B}'(A')\oplus A')
	\end{eqnarray*}
	carries an anti-flexible algebra structure, 
	where $R_B(A)\oplus A$ and 
	${R_B}'(A')\oplus A'$ are equipped the  algebra 
	structure given by Eq.~\eqref{eq:RB-subalgebraproduct}.
\end{proposition}
\begin{proof}
	Let $a,b\in A$. We have:
	\begin{eqnarray}\label{eq:equiv-a}
		&&( (R_B(a),a), (\varphi(R_B(a)),\varphi(a)))\ast 
		( (R_B(b),b), (\varphi(R_B(b)),\varphi(b)))=\cr
		&&(
		(R(a)\cdot R_B(b),a\cdot R_B(b)+R_B(a)\cdot b+
		\lambda a\cdot b),\\ &&
		(\varphi(R_B(a))\cdot'\varphi(R_B(b)),
		\varphi(a)\cdot'\varphi(R_B(b))+
		\varphi(R_B(a))\cdot' \varphi(b)+
		\lambda \varphi(a)\cdot'\varphi(b))\nonumber
		)
	\end{eqnarray}
	\begin{eqnarray}\label{eq:equiv-b}
		&&( 
		(R_B(a),a), 
		({R_B}'(\varphi(a)),\varphi(a))
		)\ast 
		(
		(R_B(b),b),
		({R_B}'(\varphi(b)),\varphi(b))
		)=\cr
		&&(
		(R_B(a)\cdot R_B(b),a\cdot R_B(b)+R_B(a)\cdot b+\lambda a\cdot b),\\&&
		({R_B}'(\varphi(a))\cdot' {R_B}'(\varphi(b)),
		\varphi(a)\cdot'{R_B}'(\varphi(b))+
		{R_B}'(\varphi(a))\cdot'\varphi(b)+
		\lambda \varphi(a) \cdot' \varphi(b))\nonumber
		)
	\end{eqnarray}
	\begin{itemize}
		\item 
		If $\varphi$ is a morphism of 
		the anti-flexible Rota-Baxter operators 
		$(A, \cdot, R_B)$ and $(A', \cdot', {R_B}')$, 
		then the rhs of Eq.~\eqref{eq:equiv-a}
		and Eq.~\eqref{eq:equiv-b} are the same and can be 
		written as
		\begin{eqnarray}\label{eq:equiv-c}
			&&((
			R_B(a\cdot R_B(b)+
			R_B(a)\cdot b+\lambda a\cdot b), 
			a\cdot R_B(b)+
			R_B(a)\cdot b+\lambda a\cdot b), \cr&&
			(
			\varphi(R_B(a\cdot R_B(b)+
			R_B(a)\cdot b+\lambda a\cdot b)),
			\varphi(a\cdot R_B(b)+
			R_B(a)\cdot b+\lambda a\cdot b))
			).
		\end{eqnarray}
		Hence $	\{
		((R_B(a), a), (\varphi(R_B(a)), \varphi(a)))|  a\in A\}$ is anti-flexible sub-algebra.
		\item
		Conversely, if   $	\{
		((R_B(a), a), (\varphi(R_B(a)), \varphi(a)))|  a\in A\}$ is anti-flexible sub-algebra,
		where $R_B(A)\oplus A$ and 
		${R_B}'(A')\oplus A'$ are equipped with semi-direct
		product algebra structure given by Eq.~\eqref{eq:RB-subalgebraproduct}, 
		then the rhs of Eq.~\eqref{eq:equiv-a} is equal to 
		Eq.~\eqref{eq:equiv-c} and also equal to 
		the rhs of Eq.~\eqref{eq:equiv-b}. All theses 
		imply that $\varphi$ is an
		anti-flexible algebra morphism and satisfy Eq.~\eqref{eq:morphism-RB}, which means 
		that $\varphi$ is a morphism of the anti-flexible
		Rota-Baxter algebras
		$(A, \cdot, R_B)$ and 
		$(A', \cdot', {R_B}')$.
	\end{itemize}
	Therefore, holds the equivalence.
\end{proof}
\begin{theorem}\label{thm:pre-anti-RB}
	Let $(A, \cdot)$ be an anti-flexible algebra
	and $R_B:A\rightarrow A$ be 
	a Rota-Baxter operator of weight $\lambda$ on $A$.
	Then $R_B$ induces a pre-anti-flexible algebra 
	structure on $A$ by
	\begin{subequations}
		\begin{eqnarray}\label{eq:Rot-pre-anti-flexible1}
			a\prec b= a\cdot R_B(b)+\frac{\lambda}{2} a\cdot b, \; \forall a, b\in A,
		\end{eqnarray}
		\begin{eqnarray}\label{eq:Rot-pre-anti-flexible2}
			a\succ b= R_B(a)\cdot b+\frac{\lambda}{2}  a\cdot b, \; \forall a, b\in A,
		\end{eqnarray}
	\end{subequations}
	if and only if 
	\begin{eqnarray}\label{eq:equiv}
		\lambda^2( (a\cdot b)\cdot c+c\cdot (b\cdot a))=0,\quad \forall a, b\in A. 
	\end{eqnarray}
\end{theorem}
\begin{proof}
	Let $a, b, c\in A$.  We have:
	\begin{eqnarray*}
		(a\succ b)\prec c-a\succ (b\prec c)
		&=&
		(R_B(a), b, R_B(c))+
		\frac{\lambda}{2}(R_B(a), b, c)\cr&+&
		\frac{\lambda}{2}(a, b, R_B(c))+
		\frac{\lambda^2}{4}(a, b,c), 
	\end{eqnarray*}
	\begin{eqnarray*}
		(a\succ b+a\prec b)\succ c-a\succ (b\succ c)
		&=&
		(R_B(a), R_B(b), c)+
		\frac{\lambda}{2}(a, R_B(b), c)+
		\frac{\lambda}{2}(R_B(a), b, c)\cr&+&
		\frac{\lambda^2}{4}(a, b, c)+
		\frac{\lambda^2}{4}(a\cdot b)\cdot c,
	\end{eqnarray*}
	\begin{eqnarray*}
		(c\prec b)\prec a-
		c\prec (b\succ a+b\prec a)
		&=&
		(c, R_B(b), R_B(a))+
		\frac{\lambda}{2}(c, b, R_B(a))+
		\frac{\lambda}{2}(c, R_B(b), a)\cr&+&
		\frac{\lambda^2}{4}(c, b, a)-
		\frac{\lambda^2}{4}c\cdot (b\cdot a).
	\end{eqnarray*}
	Since $(A, \cdot)$ is an anti-flexible algebra, then 
	$(A, \prec, \succ)$ is a pre-anti-flexible algebra if and only if 
	$\forall a, b, c\in A, \lambda^2((a\cdot b)\cdot c+c\cdot (b\cdot a))=0.$
\end{proof}
\begin{remark}
	In view of Theorem~\ref{thm:pre-anti-RB}, 
	a Rota-Baxter operator of weight zero  $R_B$  
	induces a pre-anti-flexible algebra 
	structure on $A$ which is given by
	\begin{eqnarray}
		a\prec b= a\cdot R_B(b), \;
		a\succ b= R_B(a)\cdot b, \; \forall a, b\in A.
	\end{eqnarray}
	This result was established in \cite[Theorem 2.2 \& Corollary 2.3]{Mafoya}.
\end{remark}
\begin{proposition}
	Let $(A, \cdot)$ be an anti-flexible algebra
	and $R_B:A\rightarrow A$ be 
	a Rota-Baxter operator of weight $\lambda$ on $A$.
	Then both $(A, \ast)$ and $(A, \star)$ are left symmetric,
	right symmetric algebras, respectively, where  
	$\forall a, b\in A,$
	\begin{subequations}
		\begin{eqnarray}
			a\ast b=[R_B(a), b]+\frac{\lambda}{2} [a,b],
		\end{eqnarray}
		\begin{eqnarray}
			a\star b=[a, R_B(b)]+\frac{\lambda}{2} [a,b],
		\end{eqnarray}
	\end{subequations}
	where for $a,b\in A, [a,b]=a\cdot b-b\cdot a,$
	if and only if the relation Eq.~\eqref{eq:equiv} is satisfied. 
\end{proposition}
\begin{proof}
	By straightforward calculation, we get the result.
\end{proof}
\begin{theorem}
	Let $(A, \ast)$ be an anti-flexible equipped with 
	a non degenerate skew symmetric  bilinear form $\omega: A\times A \rightarrow \mathbb{K}$ 
	satisfying, 
	\begin{eqnarray}
		\omega(a\ast b, c)+\omega (b\ast c, a)+\omega(c\ast a, b)=0, \forall a,b,c\in A.
	\end{eqnarray}
	There is a pre-Lie algebra structure $``\circ"$ on $A$ given by 
	$\forall a,b,c\in A,$
	\begin{eqnarray}\label{eq:pre-Lie}
		\omega(a\circ b, c)=\omega(b, [c,a]),
	\end{eqnarray} 
	where $[a,b]=a\ast b-b\ast a.$
\end{theorem}
\begin{proof}
	For any $a,b,c\in A,$ we have
	\begin{eqnarray*}
		\omega([a,  b], c)+\omega ([b,  c], a)+\omega( [c,  a], b)&=&
		\omega(a\ast b, c)+\omega (b\ast c, a)\cr&&+\omega(c\ast a, b)-
		\omega(b\ast a, c)\cr&&-\omega(a\ast c, b)-\omega (c\ast b, a)=0.
	\end{eqnarray*}
	Therefore, 
	$\omega$ is a non degenerate skew symmetric  bilinear on the Lie algebra $(A, [,])$
	satisfying the following identity, 
	\begin{eqnarray}
		\omega([a,  b], c)+\omega ([b,  c], a)+\omega( [c,  a], b)=0, \forall a,b,c\in A.
	\end{eqnarray}
	Thus, $(A, [,], \omega)$ is a symplectic Lie algebra, hence 
	there is a pre-Lie algebra structure ``$\circ$" on $A$ given by Eq.~\eqref{eq:pre-Lie}. 
\end{proof}
\begin{proposition}
	If $\varphi$ is a morphism 
	of the Rota-Baxter anti-flexible algebras 
	$(A, \cdot, R_B)$ and $(A', \cdot', {R_B}')$, then it 
	is a morphism of its induced pre-anti-flexible algebras. 
\end{proposition}
\begin{proof}
	For any  $a,b\in A$, we have:
	\begin{eqnarray*}
		\varphi(a\prec b)&=&		
		\varphi(a\cdot R_B(b)+\frac{\lambda}{2}a\cdot b)\cr&=&
		\varphi(a)\cdot' \varphi(R_B(b))+
		\frac{\lambda}{2}\varphi(a)\cdot'\varphi(b)=
		\varphi(a)\prec'\varphi(b).
	\end{eqnarray*}
	Similarly, we also have 
	$\varphi(a\succ b)=\varphi(a)\succ'\varphi(b).$
	Hence the result follows.
\end{proof}
For a given anti-flexible 
Rota-Baxter algebra $(A, \cdot, R_B)$, we have
$\forall a,b\in A$,
\begin{eqnarray*}
	[R_B(a), R_B(b)]&=&
	R_B(a)\cdot R_B(b)-
	R_B(b)\cdot R_B(a)\cr 
	&=& R_B(a\cdot R_B(b)+R_B(a)\cdot b+\lambda a\cdot b)\cr
	&-&
	R_B(b\cdot R_B(a)+R_B(b)\cdot a+\lambda b\cdot a)\cr
	[R_B(a), R_B(b)] &=&
	R_B([a, R_B(b)]+[R_B(a), b]+\lambda [a,b]).
\end{eqnarray*}

Hence we have the following proposition:
\begin{proposition}
	For a given anti-flexible 
	Rota-Baxter algebra \newline $(A, \cdot, R_B)$,
	setting by $\forall $ $a, b\in A$, 
	$[a, b]=a\cdot b-b\cdot a$, then
	$(A, [,], R_B)$ is a Lie 
	Rota-Baxter algebra.
\end{proposition}
\begin{proof}
	By straightforward calculation, we get the result.
\end{proof}
\begin{theorem}
	Let $(A, \cdot)$ be an anti-flexible algebra such that for 
	any $a,b,c\in A,$
	\begin{eqnarray}
		(a\cdot b)\cdot c+c\cdot (b\cdot a)=0.
	\end{eqnarray}
	Consider $R_B:A\rightarrow A$ a linear map 
	on $A.$ The linear products given by Eqs.~\eqref{eq:Rot-pre-anti-flexible1} and \eqref{eq:Rot-pre-anti-flexible2}, i.e.,
	$ \forall a, b\in A,$
	\begin{subequations}
		\begin{eqnarray*}
			a\prec b= a\cdot R_B(b)+\frac{\lambda}{2} a\cdot b, 
		\end{eqnarray*}
		\begin{eqnarray*}
			a\succ b= R_B(a)\cdot b+\frac{\lambda}{2}  a\cdot b, 
		\end{eqnarray*}
	\end{subequations}
	
	defines an anti-flexible algebra structure on $A$
	if and only if $R_B$ is a Rota-Baxter operator of weight $\lambda$ on $A.$
\end{theorem}
\begin{proof}
	Using the Theorem~\ref{thm:pre-anti-RB}, we easily establish the result.  
\end{proof}
\section{$\mathcal{O}$-operators on anti-flexible algebras}\label{section3}
In this section, we recall the notion of 
$\mathcal{O}$-operators on anti-flexible
and establish their basic properties, we 
use an $\mathcal{O}$-operator on anti-flexible
to construct  pre-anti-flexible and prove that 
any  $\mathcal{O}$-operator on anti-flexible induces
left-symmetric algebra, and related consequences are derive.
\begin{definition}
	Let $(A, \cdot)$ be an anti-flexible algebra. A vector space $M$
	is said to be an $A$-bimodule if $A$ induces two actions \newline ``left" 
	$\begin{array}{llll}
		A\times M&\rightarrow& M \\
		(a, m)&\mapsto&  am
	\end{array}
	$ and 
	``right" 
	$\begin{array}{llll}
		M\times A&\rightarrow& M \\
		(m, a)&\mapsto&  ma
	\end{array}
	$
	in $M$
	and satisfying: 
	\begin{subequations}
		\begin{eqnarray}\label{eq:bimodule1a}
			(a\cdot b)m-a(bm)-(mb)a+m(b\cdot a)=0,\; \; \forall a,b\in A,  m\in M,
		\end{eqnarray}
		\begin{eqnarray}\label{eq:bimodule2a}
			a(mb)-(am)b-b(ma)+(bm)a=0, \;\; \forall a,b\in A,  m\in M.
		\end{eqnarray}
	\end{subequations} 
\end{definition}
\begin{remark}
	For a given  bimodule $M$ over an anti-flexible algebra $(A, \cdot)$,
	there are two linear maps $l: A\times M\rightarrow M, (a, m)\mapsto l(a, m):=am$ and 
	$r: M\times A\rightarrow M, (m, a)\mapsto r(m, a):=ma$
	satisfying:  $\forall a,b\in A\mbox{ and }\forall m\in M,$
	\begin{subequations}
		\begin{eqnarray}\label{eq:bimodule1}
			l(a\cdot b, m)-l(a, l(b, m))=r(r(m, b), a)-r(m, b\cdot a), 
		\end{eqnarray}
		\begin{eqnarray}\label{eq:bimodule2}
			l(a, r(m, b))-r(l(a, m), b)=
			l(b, r(m, a))-r(l(b, m), a).
		\end{eqnarray}
		Clearly, Eq.~\eqref{eq:bimodule1a} corresponds to Eq.~\eqref{eq:bimodule1} while 
		Eq.~\eqref{eq:bimodule2a} corresponds to Eq.~\eqref{eq:bimodule2}.
	\end{subequations}
\end{remark}
\begin{example}
	The dual space $A^*$ of $A$ equipped with
	$l^*:A^*\otimes A\rightarrow A^*$ and 
	$r^*:A\otimes A^*\rightarrow A^*$ given by:
	\begin{eqnarray*}
		l^*(a,f)(b):=f(b\cdot a) \mbox{  and  } 
		r^*(f, a)(b):=f(a\cdot b), \forall a,b\in A, f\in A^*,
	\end{eqnarray*}
	carries an $A$-bimodule 
	structure, i.e.,  $\forall a, b,c\in A$ and  $\forall f\in A^*,$
	\begin{eqnarray*}
		l^*(a\cdot b, f)(c)-l^*(a, l^*(b,f))(c)&=&
		f(c\cdot(a\cdot b))-l^*(b,f)(c\cdot a)\cr&=&f(c\cdot (a\cdot b)-(c\cdot a)\cdot b)\cr 
		&=&f(b\cdot (a\cdot c)-(b\cdot a)\cdot c)\cr&=&
		r^*(r^*(f, b), a)(c)-r^*(f, b\cdot a)(c)
	\end{eqnarray*}
	Hence,
	$l^*(a\cdot b, f)-l^*(a, l^*(b,f))=r^*(r^*(f, a), b)-r^*(f, b\cdot a)$. 
	Thus Eq.~\eqref{eq:bimodule1} is satisfied.
	\begin{eqnarray*}
		l^*(a, r^*(f, b))(c)-r^*(l^*(a, f), b)(c)&=&
		r^*(f, b)(c\cdot a)-l^*(a, f)(b\cdot c)\cr&=&
		f(a\cdot (c\cdot b)-(a\cdot c)\cdot b)\cr&=&
		l^*(b, r^*(f, a))(c)-r^*(l^*(b,f), a)(c)
	\end{eqnarray*}
	Hence, 
	$l^*(a, r^*(f, b))-r^*(l^*(a, f), b)=l^*(b, r^*(f, a))-r^*(l(b,f), a)$.
	Thus Eq.~\eqref{eq:bimodule2} is satisfied. 
\end{example} 
\begin{proposition}\label{Prop-Sum}
	For a given anti-flexible algebra $(A, \cdot)$ 
	and an $A$-bimodule $M,$ there is an anti-flexible algebra 
	structure ``$\ast$" on the vector space $A\oplus M$ given by 
	\begin{eqnarray}\label{eq:subalgebraproduct}
		(a,m)\ast (b, n)=
		(a\cdot b, l(a,n)+r(m, b)), \forall a, b\in A, \forall m, n \in M.
	\end{eqnarray}
\end{proposition}
\begin{proof}
	Let $a, b, c\in A$ and $m, n, s\in M$. The associator of 
	the bilinear product given by Eq.~\eqref{eq:subalgebraproduct}
	is derived by:
	\begin{eqnarray*}
		[ (a,m), (b, n), (c, s)]_{_\ast}&=&
		([a,b,c], l(a\cdot b, s)-l(a, l(b,s))
		+r(l(a,n), c)\cr&-&l(a, r(n, c))+r(r(m,b), c)-r(m, b\cdot c))
	\end{eqnarray*}
	and similarly
	\begin{eqnarray*}
		[ (c,s), (b, n), (a, m)]_{_\ast}&=&
		([c,b,a], l(c\cdot b, m)-l(c, l(b,m))
		+r(l(c,n), a)-l(c, r(n, a))\cr&&+r(r(s,b), a)-r(s, b\cdot a)).
	\end{eqnarray*}
	According to Eq.~\eqref{eq:bimodule1} 
	and Eq.~\eqref{eq:bimodule2}, we deduce that 
	$\forall a, b, c\in A$ and $m, n, s\in M$,
	$[ (a,m), (b, n), (c, s)]_{_\ast}=[ (c,s), (b, n), (a, m)]_{_\ast}$, 
	hence $(A\oplus M, \ast)$ is an anti-flexible algebra. 
\end{proof}
\begin{definition}
	Let $(A, \cdot)$ be an anti-flexible algebra and $M$ be an $A$-bimodule. 
	A linear map $T:M\rightarrow A$ is said to be an $\mathcal{O}$-operator 
	on $A$ with respect to the bimodule $M$ if  
	\begin{eqnarray}\label{eq:o-operator}
		T(m)\cdot T(n)=T(mT(n)+T(m)n),\forall m, n\in M.
	\end{eqnarray}  
\end{definition}
\begin{proposition}
	A linear map $T: M\rightarrow A$ is an $\mathcal{O}$-operator 
	on $A$ with respect to the $A$-bimodule $M$ if and only if
	the bilinear product given by
	\begin{eqnarray}\label{eq:anti-flex-o-operator}
		m\star n=mT(n)+T(m)n, \quad \forall m, n\in M,
	\end{eqnarray}
	defines an anti-flexible algebra structure on $M.$ 
	Furthermore,  the anti-flexible algebras
	$(M, \star)$ and $(T(M), \cdot)$ are homomorphic.
\end{proposition}
\begin{proof}
	Let $m,n,p\in M.$ Using the equation Eq.~\eqref{eq:o-operator}, 
	then the associator of related to the bilinear product given in Eq.~\eqref{eq:anti-flex-o-operator}
	can been written in the following form
	\begin{eqnarray*}
		[m, n, p]_{\star} &=&(m\star n)\star p-m\star(n\star p)
		=
		(mT(n)+T(m)n)T(p)\cr&&
		+T(mT(n)+T(m)n)p
		-mT(nT(p)+T(n)p)- T(m)(nT(p)+T(n)p)\cr
		&=&(mT(n))T(p)
		+(T(m)n)T(p)
		+(T(m)\cdot T(n))p\cr
		&&-m(T(n)\cdot T(p))
		-T(m)(nT(p))
		-T(m)(T(n)p)\cr
		&=&
		\{(T(m)\cdot T(n))p
		-T(m)(T(n)p)\}
		+\{(mT(n))T(p)
		-m(T(n)\cdot T(p))\}
		\cr&&+
		\{(T(m)n)T(p)
		-T(m)(nT(p))\}
	\end{eqnarray*}
	Thus we have
	\begin{eqnarray}
		[m, n, p]_{\star}&=&
		\{(T(m)\cdot T(n))p-T(m)(T(n)p)\}
		+\{(mT(n))T(p)-m(T(n)\cdot T(p))\}
		\\ &&+
		\{(T(m)n)T(p)-T(m)(nT(p))\}\nonumber
	\end{eqnarray}
	Furthermore, $\forall m, n, p \in M,$ we have
	\begin{eqnarray*}
		[m,n, p]_{\star}-[p, n, m]_{\star} &=&
		\{(T(m)\cdot T(n))p-T(m)(T(n)p)-
		(pT(n))T(m)+p(T(n)\cdot T(m))\}\cr &&+
		\{(mT(n))T(p)-m(T(n)\cdot T(p))-
		(T(p)\cdot T(n))m+T(p)(T(n)m)\}\cr&&+
		\{(T(m)n)T(p)-T(m)(nT(p))-
		(T(p)n)T(m)+T(p)(nT(m))\}.
	\end{eqnarray*}
	Using the equations Eq.~\eqref{eq:bimodule1a} and Eq.~\eqref{eq:bimodule2a}
	on the right hand sid of the above identity, we then get
	\begin{equation}
		(m\star n)\star p-m\star(n\star p)=
		(p\star n)\star m-p\star(n\star m), \forall m, n, p \in M,
	\end{equation}
	which proves that $(M, \star)$ is an anti-flexible algebra.
	Besides, we also have:
	$T(m\star n) =T(m)\cdot T(n), \forall m, n\in M,$ thus
	$(M, \star)$ is homomorphic to $(T(M), \cdot)$
\end{proof}
\begin{proposition}
	A linear map $T: M\rightarrow A$ is an $\mathcal{O}$-operator 
	on $A$ with respect to the $A$-bimodule $M$ if and only if the set  
	$
	\{(T(m), m)| m\in M\}
	$
	equipped with the semi-direct product given by Eq.~\eqref{eq:subalgebraproduct} 
	is a sub-algebra of the semi-direct product algebra $A\oplus M$.
\end{proposition}
\begin{proof}
	By the  definition of $T$, we have  $T(M)\subseteq A$, thus 
	$
	\{(T(m), m)| m\in M\}
	\subseteq A\oplus M.$
	
	For any $m, n \in M, $
	\begin{eqnarray*}
		(T(m), m)\ast (T(n), n)=
		(T(m)\cdot T(n), mT(n)+T(m)n).
	\end{eqnarray*}
	According to the definition of 
	semi-direct production given by Eq.~\eqref{eq:subalgebraproduct}, 
	the above bilinear product ``$\ast$" 
	is stable on
	$
	\{(T(m), m)| m\in M\}
	$ if and only if 
	$\forall m, n\in M,$
	\begin{eqnarray*}
		T(m)\cdot T(n)=T(T(m)n+mT(n)),
	\end{eqnarray*} hence 
	$T$ is an $\mathcal{O}$-operator 
	on $A$ with respect to 
	the $A$-bimodule $M$.
\end{proof}
\begin{proposition}\label{prop:pre}
	Let $T: M\rightarrow A$ be an 
	$\mathcal{O}$-operator on $A$ with 
	respect to the bimodule $M$.
	There is a pre-anti-flexible algebra
	structure on $M$ given by $ \forall m,n\in M$
	\begin{eqnarray}\label{eq:product-pre-anti-flexible}
		m\prec n=r(m, T(n))=mT(n), \quad m\succ n=l(T(m), n)=T(m)n.
	\end{eqnarray}
\end{proposition}
\begin{proof}
	Let $m, n, s \in M.$ Using 
	Eq.~\eqref{eq:product-pre-anti-flexible}, 
	Eq.~\eqref{eq:bimodule1}, Eq.~\eqref{eq:bimodule2},
	and Eq.~\eqref{eq:o-operator}, we have
	\begin{eqnarray*}
		(m\succ n)\prec s-m\succ (n\prec s)&=&
		(T(m)n)T(n)-T(m)(nT(s))\cr&=&
		r(l(T(m),n), T(s))-l(T(m), r(n, T(s)))\cr&=&
		r(l(T(s),n), T(m))-l(T(s), r(n, T(m)))\cr&=&
		(s\succ n)\prec m-s\succ (n\prec m).
	\end{eqnarray*}
	Besides, we have 
	\begin{eqnarray*}
		(m\prec n +m\succ n)\succ s-m\succ (n\succ s)&=&
		T(mT(n)+T(m)n)s -T(m)(T(n)s)\cr&=&
		(T(m)T(n))s\ -T(m)(T(n)s)\cr&=&
		l(T(m)T(n))s -l(T(m))(l(T(n),s))\cr&=&
		r(r(s,T(n)), T(m)) -r(s, T(n)T(m))\cr&=&
		(sT(n))T(m) -s(T(nT(m)+T(n)m))\cr&=&
		(s\prec n)\prec m-
		s\prec(n\prec m+n\succ m).
	\end{eqnarray*} 
	Hence  $(M, \prec, \succ)$ is a
	pre-anti-flexible algebra.
\end{proof}
\begin{theorem}
	Let $T: M\rightarrow A$ be an 
	$\mathcal{O}$-operator on an anti-flexible algebra $(A, \cdot)$ with 
	respect to the bimodule $M$. Then $T$ induces both left-symmetric and right-symmetric
	algebras structures on $M$, respectively by  $\forall n, m\in M, $
	\begin{eqnarray*}
		m\star n=T(m)n-nT(m),
	\end{eqnarray*}
	\begin{eqnarray*}
		m\ast n=nT(m)-T(m)n.
	\end{eqnarray*}
\end{theorem}
\begin{proof}
	By considering  Proposition~\ref{prop:pre} and  Theorem~\ref{thm:Prelie}, then follows the result.
\end{proof}
\begin{definition}
	A morphism of two $\mathcal{O}$-operators $T$ and $T'$
	with respect to the $A$-bimodule 
	$M$ of anti-flexible algebra $(A, \cdot_{A})$
	and to the $B$-bimodule 
	$N$ of anti-flexible algebra $(B, \cdot_{B})$, 
	respectively, is a pair $(\phi, \varphi)$ of 
	an anti-algebra morphism $\phi: A\rightarrow B$
	and a linear map $\varphi: M\rightarrow N$
	such that the following diagrams commute
	\begin{subequations}
		\begin{eqnarray}\label{eq:isomorphism-cond1}
			\begin{array}{cccc}
				\xymatrix{	
					& M\ar[rr]^-{T}\ar[d]^-{\varphi}
					&& A \ar[d]^-{\phi}\cr
					& N\ar[rr]^-{T'}
					&& B
				}
			\end{array}
			\mbox{  equivalently }  
			\phi \circ T=T'\circ \varphi,
		\end{eqnarray}
		\begin{eqnarray}\label{eq:isomorphism-cond2}
			\begin{array}{cccc}
				\xymatrix{	
					& A\oplus M \ar[d]^-{\phi\oplus\varphi}\ar[rr]^-{l}
					&& M \ar[d]^-{\varphi}\\
					& B\oplus N\ar[rr]^-{l'}&& N
				}
			\end{array}
			\mbox{ equivalently } \cr
			\varphi(l(a, m))=l'(\phi(a), \varphi(m)), 
		\end{eqnarray}
		\begin{eqnarray}\label{eq:isomorphism-cond3}
			\begin{array}{cccc}
				\xymatrix{	
					& M\oplus A \ar[d]^-{\varphi\oplus\phi}\ar[rr]^-{r}
					&& M \ar[d]^-{\varphi}\\
					& N\oplus B\ar[rr]^-{r'}&& N
				}
			\end{array}
			\mbox{    equivalently }\cr 
			\varphi(r(m ,a))=
			r'(\varphi(m), \phi(a)).
		\end{eqnarray}
	\end{subequations}
	Particularly, if both 
	$\phi$ and $\varphi$
	are linear isomorphisms, 
	then $(\phi, \varphi)$
	is an isomorphism of the 
	$\mathcal{O}$-operators $T$
	and $T'$.
\end{definition}
\begin{proposition}
	Let $(A, \cdot_1)$ and $(B, \cdot_2)$ be two anti-flexible algebras, $M$ be a $A$-bimodule, 
	and $N$ be a $B$-bimodule.	
	Consider	
	$(\phi, \varphi)$ be a pair of linear maps 
	$\phi: A\rightarrow B$ and  
	$\varphi: M\rightarrow N,$ 
	$T: M\rightarrow A$ be an $\mathcal{O}$-operator on 
	$A$ with respect to $M$, and $T': N\rightarrow B$ be an $\mathcal{O}$-operator on 
	$B$ with respect to $N.$
	The following conditions are equivalent: 
	\begin{enumerate}
		\item 
		The pair $(\phi, \varphi)$ is a morphism of the 
		$\mathcal{O}$-operators $T$ to $T'.$ 
		\item 
		The set
		\begin{eqnarray}\label{eq:set2}
			\{
			((T(m), m), (\phi(T(m)), \varphi(m))
			)|  m\in M
			\}	\subset (A\oplus M)\oplus (B\oplus N)
		\end{eqnarray}
		is an anti-flexible sub-algebra, where $A\oplus M$ and 
		$B\oplus N$ are equipped with semi-direct
		product algebra structure.
	\end{enumerate}
\end{proposition}
\begin{proof}
	Let $m, n\in M$, we have:
	\begin{eqnarray}\label{eq:equiv1}
		&&((T(m), m), (\phi(T(m)), \varphi(m)))\ast  
		((T(n), n), (\phi(T(n)), \varphi(n)))=\cr&&
		( (T(m)\cdot_1T(n), T(m)n+mT(n)), \cr&&
		(    \phi(T(m))\cdot_2\phi(T(n)),\phi(T(m))\varphi(n)+ \varphi(m)\phi(T(n))  )
		)
	\end{eqnarray}
	and 
	\begin{eqnarray}\label{eq:equiv2}
		&&(
		(T(m), m), (T'(\varphi(m)), \varphi(m))
		)\ast
		(
		(T(n), n), (T'(\varphi(n)), \varphi(n))
		)=
		\cr &&
		( (T(m)\cdot_1T(n), T(m)n+mT(n)), \cr&&
		( T'(\varphi(m))\cdot_2  T'(\varphi(n)),
		T'(\varphi(m))\varphi(n)+\varphi(m)T'(\varphi(n))
		).
	\end{eqnarray}
	\begin{enumerate}
		\item
		If the pair $(\phi, \varphi)$ of linear maps 
		$\phi: A\rightarrow B$ and  
		$\varphi: M\rightarrow N$
		is a morphism of the
		$\mathcal{O}$-operators $T$ and $T'$, 
		then we easily get the rhs of Eq.~\eqref{eq:equiv1}
		is equal to the rhs of Eq.~\eqref{eq:equiv2}
		and can be expressed as 
		\begin{eqnarray*}
			((T( T(m)n+mT(n)),  T(m)n+mT(n)),\cr
			(\phi(T(T(m)n+mT(n))), \varphi(T(m)n+mT(n)))).
		\end{eqnarray*}
		Hence, the set given in Eq.~\eqref{eq:set2}
		is a sub-algebra of 
		$(A\oplus M)\oplus (B\oplus N)$, where $A\oplus M$ and 
		$B\oplus N$ are equipped with the  semi-direct
		product algebra structure, respectively.
		\item 
		Conversely, if the set given in Eq.~\eqref{eq:set2} is a sub-algebra of 
		$(A\oplus M)\oplus (B\oplus N)$, where $A\oplus M$ and 
		$B\oplus N$ are respectively equipped with the  semi-direct
		product algebra structure, then using 
		the equations Eq.~\eqref{eq:equiv1} and  Eq.~\eqref{eq:equiv2}, the fact that 
		$T$ and $T'$ are $\mathcal{O}$-operators, we establish by direct calculation
		that $\phi: A\rightarrow B$ is an algebra morphism such that 
		both with the linear map $\varphi: M\rightarrow N$ satisfy
		Eq.~\eqref{eq:isomorphism-cond1}, Eq.~\eqref{eq:isomorphism-cond2} 
		and Eq.~\eqref{eq:isomorphism-cond3}. 
		Consequently,  $(\phi, \varphi)$ of linear maps 
		$\phi: A\rightarrow B$ and  
		$\varphi: M\rightarrow N$
		is a morphism of the
		$\mathcal{O}$-operators $T$ and $T'$. 
	\end{enumerate}
	Hence, holds the equivalence.
\end{proof}
\begin{proposition}
	Let $T:M\rightarrow A$ be an $\mathcal{O}$-operator  on $A$
	with respect to a bimodule $M$ and $T'$
	be an $\mathcal{O}$-operator on $B$ with respect 
	to a bimodule $N$.
	If $(\phi, \varphi)$ is a morphism from 
	$T$ to $T'$, then  $\varphi: M \rightarrow N$
	is a morphism between induced 
	pre-anti-flexible algebra structures derived on $M$ and $N$  by 
	Eq.~\eqref{eq:product-pre-anti-flexible}.
\end{proposition}
\begin{proof}
	For any $m, n\in M$, we have
	\begin{eqnarray*}
		\varphi(m\prec_{_M} n)
		=
		\varphi(r(m, T(n)))
		=
		r(\varphi(m), \phi(T(n)))
		=
		r(\varphi(m), T'(\varphi(n)))
		=\varphi(m)\prec_{_N}\varphi(n),
	\end{eqnarray*}
	\begin{eqnarray*}
		\varphi(m\succ_{_M} n)
	=
		\varphi(l(T(m), n))
		=
		l(\phi(T(m)), \varphi(n))
		=
		l(T'(\varphi(m)), \varphi(n))
		=\varphi(m)\succ_{_N}\varphi(n).
	\end{eqnarray*}
	Hence $\varphi:M\rightarrow N$ is a 
	morphism between
	the pre-anti-flexible algebras $(M, \prec_{_M}, \succ_{_M})$
	and $(N, \prec_{_N}, \succ_{_N})$, where  
	\begin{eqnarray*}
		&&m\prec_{_M} n=r(m, T(n)), \quad 
		m\succ_{_M} n= l(T(m), n), \forall m, n\in M, \\
		&&m'\prec_{_N}n'=r(m', T'(n')),\quad 
		m'\succ_{_N} n'=l(T'(m'), n'), \forall m', n'\in N.
	\end{eqnarray*} 
\end{proof}
\begin{lemma}
	Let $(A, \cdot)$ be an anti-flexible algebra, 
	$M$ be a bimodule of $A$ and $T: M\rightarrow A$
	an $\mathcal{O}$-operator on $A$
	with respect to the bimodule $M$.
	The following linear maps 
	\begin{subequations}
		\begin{eqnarray}
			l_{_T}: M\times A\rightarrow A, 
			l_{_T}(m, a)= T(m)\cdot a-T(ma),
		\end{eqnarray}
		\begin{eqnarray}
			r_{_T}: A\times M\rightarrow A, 
			r_{_T}(a, m)=a\cdot  T(m)-T(am), 
		\end{eqnarray}
	\end{subequations}
	$\forall a\in A, m\in M,$	define a $M$-bimodule structure on $A$.
\end{lemma}
\begin{proof}
	Let $m, n\in M$ and $a\in A$.
	Since $T$ is 		an $\mathcal{O}$-operator on $A$
	with respect to the bimodule $M$, we have:	
	\begin{eqnarray*}
		l_{_T}(m\star n, a)-l_{_T}(m, l_{_T}(n, a))&=&
		(T(m)\cdot T(n))\cdot a-T((mT(n)+T(m)n)a)\cr&-&T(mT(na))
		+T(m(T(n)\cdot a))\cr&+& T(m)\cdot T(na)-T(m)\cdot(T(n)\cdot a)\cr&=&
		(T(m), T(n), a)-T((mT(n))a)\cr&-&T((T(m)n)a)
		+\sout{T(mT(na))}\cr&+& T(T(m)(na))- \sout{T(mT(na))}+ T(m(T(n)\cdot a))\cr&=&
		(T(m), T(n), a)-T((mT(n))a)\cr&-& T((T(m)n)a)+
		T(T(m)(na))+ T(m(T(n)\cdot a))
	\end{eqnarray*}
	and 
	\begin{eqnarray*}
		r_{_T}(r_{_T}(a,n), m)-r_{_T}(a, n\star m)&=&
		(a\cdot T(n))\cdot T(m)-a\cdot(T(n)\cdot T(n))\cr&-& T(an)\cdot T(m)+
		T(T(an)m)\cr&+& T(a(nT(m)))+T(a(T(n)m))- T((a\cdot T(n))m)\cr&=&
		(a, T(n), T(m))-T((a\cdot T(n))m)\cr&-& T((an)T(m))-
		T(T(an)m)\cr&+& T(T(an)m)+T(a(nT(m)))+ T(a(T(n)m))\cr&=&
		(a, T(n), T(m))-T((a\cdot T(n))m)\cr&-& T((an)T(m))+
		T(a(nT(m)))+ T(a(T(n)m))
	\end{eqnarray*}
	In view of $(A, \cdot)$ is an anti-flexible algebra and $M$ is an 
	$A$-module, we have:
	\begin{eqnarray*}
		&&l_{_T}(m\star n, a)-l_{_T}(m, l_{_T}(n, a))
		-r_{_T}(r_{_T}(a,n), m)+r_{_T}(a, n\star m)=\cr
		&&T(m(T(n)\cdot a))
		-T((mT(n))a)
		-T(a(T(n)m))
		+T((a\cdot T(n))m)\cr&& +
		T(T(m)(na))
		-T((T(m)n)a) 
		+T((an)T(m))
		-T(a(nT(m)))=\cr&&
		T(
		r(m, T(n)\cdot a)
		-r(r(m, T(n)), a)
		+l(a\cdot T(n), m)
		\cr&&-l(a, l(T(n), m))
		)+
		T(
		l(T(m), r(n, a))
		-r(l(T(m), n), a)
		\cr&&+r(l(a, n), T(m))
		-l(a, r(n, T(m)))
		)=0.
	\end{eqnarray*}
	Thus,  $\forall a\in A$ and  $\forall m, n\in M$
	\begin{eqnarray}\label{eq:extend-bimodule1}
		l_{_T}(m\star n, a)-l_{_T}(m, l_{_T}(n, a))=
		r_{_T}(r_{_T}(a,n), m)-r_{_T}(a, n\star m).
	\end{eqnarray}
	Similarly, we have  $\forall a\in A$ and
	$\forall m, n\in M$, 
	\begin{eqnarray}\label{eq:extend-bimodule2}
		l_{_T}(m,r_{_T}(a, n))-	
		r_{_T}(l_{_T}(m, a), n)	
		=	
		l_{_T}(n,r_{_T}(a, m))-	
		r_{_T}(l_{_T}(n, a), m).	
	\end{eqnarray}
	Hence $(l_{_T}, r_{_T})$ satisfies 
	Eq.~\eqref{eq:bimodule1} and Eq.~\eqref{eq:bimodule2}.
	Therefore  the proof.
\end{proof}
\section{Anti-flexible Nijenhuis algebras}\label{section4}
This section contains basic definitions and properties of 
anti-flexible Nijenhuis algebras.  
\begin{definition}
	Let $(A, \cdot)$ be an anti-flexible algebra. A linear map 
	$N_J:A\rightarrow A$ is said to be a Nijenhuis operator 
	if 
	\begin{eqnarray}\label{eq:Nijenhuis}
		N_J(x)\cdot N_J(y)-
		N_J(N_J(x)\cdot y+ x\cdot N_J(y)
		-N_J(x\cdot y)), \forall x, y\in A.
	\end{eqnarray}  
\end{definition}
\begin{proposition}\label{prop:Nj-stability}
	Let $(A, \cdot)$ be an anti-flexible algebra 
	and $N_J:A\rightarrow A$ be a Nijenhuis operator on $A$.
	The bilinear operation $\cdot_{_{N_J}}: A\times A \rightarrow A$ given by
	\begin{eqnarray}\label{eq:Nijenhuis-product}
		x\cdot_{_{N_J}} y=
		N_J(x)\cdot y+
		x\cdot N_J(y)-N_J(x\cdot y), \forall x, y\in A,
	\end{eqnarray}
	is such that $(A, \cdot_{_{N_J}})$ 
	is an anti-flexible algebra.
\end{proposition}
\begin{proof}
	For any $x,y,z\in A$, the associator of 
	the bilinear product ``$\cdot_{_{N_J}}$"
	can be expressed in term of the associator of 
	the bilinear product ``$\cdot$" as,  $\forall x,y,z\in A$,
	\begin{eqnarray}\label{eq:expr-assosicator-Nj}
		[x,y,z]{_{_{N_J}}}&=&
		[N_J(x), N_J(y), z]+
		[N_J(x), y, N_J(z)]+
		[x, N_J(y), N_J(z)]+
		N_J^2\left([x,y,z]\right)\cr&&-
		N_J\left([N_J(x), y, z]\right)-
		N_J\left([x, N_J(y), z]\right)-
		N_J\left([x, y, N_J(z)]\right).
	\end{eqnarray} 
	Since $(A,\cdot)$ is an
	anti-flexible algebra,
	we establish from 
	Eq.~\eqref{eq:expr-assosicator-Nj}
	that  $\forall x, y, z\in A$, 
	$[x,y,z]{_{_{N_J}}}=[z,y,x]{_{_{N_J}}}.$ 
	Hence, $(A, \cdot_{_{N_J}})$ is an anti-flexible algebra.
\end{proof}
\begin{remark}
	Let $(A, \cdot)$ be an anti-flexible algebra 
	and $N_J:A\rightarrow A$ be a Nijenhuis operator on $A$.
	\begin{itemize}
		\item 
		In view of Proposition~\ref{prop:Nj-stability}, 
		we deduce that the Nijenhuis operator $N_J$ 
		is an anti-flexible algebra 
		homomorphism between $(A, \cdot_{_{N_J}})$ and 
		$(A, \cdot)$.
		\item 
		According to Eq.~\eqref{eq:expr-assosicator-Nj}, 
		it is easy to see that if $(A, \cdot)$ is a pre-Lie 
		algebra, i.e.,  $\forall x,y,z\in A, [x,y,z]=[y, x, z]$, 
		then $(A, \cdot_{_{N_J}})$ is also a 
		pre-Lie algebra and in this case, $N_J:A\rightarrow A$
		satisfying Eq.~\eqref{eq:Nijenhuis} is a pre-Lie algebra morphism from 
		$(A, \cdot_{_{N_J}})$ to $(A, \cdot)$, this result and many other
		were derived in \cite[Propositions 4.6 \& 4.7]{Wang_S_B_L}. 
	\end{itemize}
\end{remark}
\begin{lemma}
	Let $(A, \cdot)$ be an anti-flexible algebra 
	and $N_J:A\rightarrow A$ be a Nijenhuis operator on $A$. 
	For any $p,q\in \mathbb{N}$ we have:
	\begin{itemize}
		\item 
		$(A, \cdot_{_{{N_J}^p}})$ is an anti-flexible algebra;
		\item 
		${N_{J}}^q$ is also a Nijenhuis operator on the anti-flexible algebra $(A, \cdot_{_{{N_J}^p}});$
		\item 
		The anti-flexible algebras $(A, (\cdot_{_{{N_J}^p}}){_{_{{N_J}^q}})}$ and 
		$(A, \cdot_{_{{N_J}^{p+q}}})$ coincide;
		\item 
		The anti-flexible algebras $(A, \cdot_{_{{N_J}^p}})$ and $(A, \cdot_{_{{N_J}^q}})$ are compatible, 
		that is, any linear combination of 
		``$\cdot_{_{{N_J}^p}}$" and ``$\cdot_{_{{N_J}^q}}$" still makes $A$ into an 
		anti-flexible algebra;
		\item 
		${N_{J}}^q$ is also a homomorphism from the anti-flexible algebra\\ 
		$(A, \cdot_{_{{N_J}^{p+q}}})$ to $(A, \cdot_{_{{N_J}^p}})$. 
	\end{itemize}
\end{lemma}
\begin{definition}
	A complex structure on a real Lie algebra 
	$(\mathfrak{g},[,])$ is a linear map
	$J: \mathfrak{g} \rightarrow \mathfrak{g}$
	satisfying $J^2=-{id}$ and 
	\begin{eqnarray}
		J([x,y])=[J(x), y]+[x, J(y)]+J([J(x), J(y)]),\forall x,y\in A.
	\end{eqnarray}
\end{definition}
\begin{proposition}
	Let $(A, \cdot)$ be an anti-flexible algebra.
	The linear map $J: A\oplus A \rightarrow A\oplus A$
	given by $J((x,a),(y,b))=((-y,-b),(x,a)),$ $\forall (x,a),(y,b) \in A\oplus A$
	is a complex structure on the Lie algebra
	structure derived on $A\oplus A$ by 
	\begin{eqnarray}\label{equ-a}
		[(x,a), (y,b)]=([x,y]_{_\cdot}, (L-R)(x)b-(L-R)(y)a), \forall x,y,a,b\in A,
	\end{eqnarray} 
	where $\forall x,y, a,b\in A, L(x)a= x\cdot a,$ and 
	$R(x)(a)=a\cdot x.$
\end{proposition}
\begin{theorem}
	Let $(A, \cdot)$ be an anti-flexible algebra 
	and $N_J:A\rightarrow A$ be a Nijenhuis operator on $A$.
	Then $N_J$ induces a pre-anti-flexible algebra structure
	on $A$ by,  $\forall x, y\in A,$
	\begin{subequations}
		\begin{eqnarray}\label{eq:ee1}
			x\succ y=N_J(x)\cdot y-\frac{1}{2}N_J(x\cdot y),
		\end{eqnarray}
		\begin{eqnarray}\label{eq:ee2}
			x\prec y=x\cdot N_J(y)-\frac{1}{2}N_J(x\cdot y),
		\end{eqnarray}
	\end{subequations} 
	if and only if  $\forall x,y,z\in A,$
	\begin{eqnarray}\label{eq:preLiecond}
		\frac{3}{2}N_J(N_J(z\cdot y)\cdot x+x\cdot N_J(y\cdot z))=
		N_J^2((z\cdot y)\cdot x+x\cdot (y\cdot z)).
	\end{eqnarray}
\end{theorem}
\begin{proof}
	Let $x,y,z\in A.$ We have
	\begin{eqnarray*}
		(x\succ y)\prec z-x\succ (y\prec z)-(z\succ y)\prec x+z\succ (y\prec x)=\cr
		\frac{1}{4}N_J(N_J(x\cdot y)\cdot z+z\cdot N_J(y\cdot x))
		-\frac{1}{4}N_J(N_J(z\cdot y)\cdot x+x\cdot N_J(y\cdot z))
	\end{eqnarray*}
	and 
	\begin{eqnarray*}
	(x\succ y+x\prec y)\succ z-x\succ (y\succ z)-(z\prec y)\prec x+z\prec (y\prec x+y\succ x)= \frac{1}{2}N_J(N_J(x\cdot y)\cdot z+\cr z\cdot N_J(y\cdot x))
		+\frac{1}{4}N_J(N_J(z\cdot y)\cdot x+x\cdot N_J(y\cdot z))
		-\frac{1}{2}N_J^2((z\cdot y)\cdot x+x\cdot (y\cdot z)).
	\end{eqnarray*}
	Therefore, Eq.~\eqref{eq:ee1} and Eq.~\eqref{eq:ee2} are satisfied if and 
	only if Eq.~\eqref{eq:preLiecond} is satisfied.
\end{proof}
\begin{proposition}
	Let $(A, \cdot)$ be an anti-flexible algebra 
	and $N_J:A\rightarrow A$ be a Nijenhuis operator on $A$.
	Then $A$ carries a left-symmetric algebra structure ``$\circ$" of the form, 
	$\forall x,y\in A,$ 
	\begin{eqnarray}
		x\circ y=[N_J(x),y]-\frac{1}{2} N_J([x,y]),
	\end{eqnarray} 
	where  $\forall x,y\in A, [x,y]=x\cdot y-y\cdot x,$
	if and only if Eq.~\eqref{eq:preLiecond} is satisfied.
\end{proposition}
By a direct calculation, we arrive the following results:
\begin{proposition}
	Let $(A, \cdot)$ be an anti-flexible algebra and 
	$N_J: A\rightarrow A$ a linear map.
	\begin{enumerate}
		\item If $N{_J}^2=0$, then $N_J$ is a Nijenhuis operator if and only if 
		$N_J$ is a Rota-Baxter operator of weight zero on $A$. 
		\item 
		If $N{_J}^2=N_J$, then $N_J$ is a Nijenhuis operator if and only if 
		$N_J$ is a Rota-Baxter operator of weight $-1$ on $A$. 
		\item 
		If $N{_J}^2=Id$, then $N$ is a Nijenhuis operator if and only if 
		$N_J+Id$ (resp. $N_J-Id$) is a Rota-Baxter operator of weight $-2$ (resp. $+2$) on $A$. 
	\end{enumerate}
\end{proposition}

\begin{proposition}
	Let $(A, \cdot)$ be an anti-flexible algebra, 
	$M$ is a vector space and $T: M\rightarrow A$  a linear map.
	The linear map $T$ is an $\mathcal{O}$-operator 
	on $A$ with respect to the bimodule $M$ if and only if, 
	$\forall \lambda \in \mathbb{K},$
	the linear map $R_T=
	\left(
	\begin{matrix}
		0 & T\cr 
		0 & -\lambda Id
	\end{matrix}
	\right):
	A\oplus M\rightarrow A\oplus M$
	is a Rota-Baxter operator of weight $\lambda$ on 
	the anti-flexible algebra $A\oplus M$ given in 
	Proposition~\ref{Prop-Sum}.
\end{proposition}
\begin{proposition}
	Let $(A, \cdot)$ be an anti-flexible algebra, 
	$M$ is a vector space
	and $T: M\rightarrow A$  a linear map.
	The following conditions are equivalent:
	\begin{enumerate}
		\item 
		The linear map $T$ is an $\mathcal{O}$-operator 
		on $A$ with respect to the bimodule $M.$
		
		\item 
		The linear map
		$N_T=
		\left(
		\begin{matrix}
			0 & T\cr 
			0 & 0
		\end{matrix}
		\right):
		A\oplus M\rightarrow A\oplus M$
		is a Nijenhuis operator on the semi-direct anti-flexible algebra 
		$A\oplus M$ given in Proposition~\ref{Prop-Sum}.
		
		\item 
		The linear map
		$\mathcal{N_T}=
		\left(
		\begin{matrix}
			0 & T\cr 
			0 & Id
		\end{matrix}
		\right):
		A\oplus M\rightarrow A\oplus M$
		is a Nijenhuis operator on the semi-direct anti-flexible algebra 
		$A\oplus M$ given in Proposition~\ref{Prop-Sum}.
		
	\end{enumerate}
\end{proposition}

\begin{center}
	\textbf{Acknowledgments}
\end{center}
The author thanks 
EPSRC GCRF project EP/T001968/1, 
Capacity building in Africa via technology-driven 
research in algebraic and arithmetic 
geometry (part of the Abram Gannibal Project)
for supporting his visit to 
Loughborough University during 
which a part of this paper 
was finalized.
He  also thanks all members of the Department of Mathematics 
of Loughborough University
for the warm welcome during his stay.
He would like to thank the anonymous reviewers for their valuable  comments  and suggestions.


\bigskip
\bigskip


\end{document}